\documentclass[12pt,a4paper]{article}
\usepackage{float}
\usepackage{amsmath}
\usepackage{amsfonts}
\usepackage{amssymb}
\usepackage{graphicx}%
\renewcommand{\baselinestretch}{1}
\numberwithin{equation}{section}
\newtheorem{theorem}{Theorem}[section]

\newtheorem{remark}{Remark}[section]

\begin{document}

\title{Recent advances on similarity solutions arising during free convection}
\author{BERNARD BRIGHI$\dag$  and JEAN-DAVID HOERNEL$\ddagger$ }
\date{}
\maketitle

\begin{center}
Universit\'e de Haute-Alsace, Laboratoire de Math\'ematiques et Applications
\vskip 0,1cm
4 rue des fr\`eres Lumi\`ere, 68093 MULHOUSE (France)
\vskip 1,5cm
\end{center}

\begin{abstract}
This paper reviews results about free convection near a vertical flat plate embedded in some saturated porous medium. We focus on a third order autonomous differential equation that gives a special class of solutions called similarity solutions. Two cases are under consideration: in the first one we prescribe the temperature on the plate and in the second one we prescribe the heat flux on it. We will also see that the same equation appears in 
other industrial processes.
\end{abstract}

\renewcommand{\baselinestretch}{1}
\footnotetext{\hspace{-0.8cm} AMS 2000 Subject Classification: 34B15, 34C11, 76D10.}
\footnotetext{\hspace{-0.8cm} Key words and phrases: Boundary layer, similarity solution, third order nonlinear differential equation,  boundary value problem. }
\footnotetext{\hspace{-0.8cm} $\dag$ bernard.brighi@uha.fr $\ddagger$ j-d.hoernel@wanadoo.fr}

\section{Introduction}
Free convection boundary layer flows near a vertical flat plate embedded in some porous medium are studied for many years and a natural way to describe the convective flow is to look for similarity solutions. We consider two different sets of boundary conditions for the temperature on the plate: either we prescribe the temperature or we prescribe the heat flux. Both cases are leading to the same following third order non-linear autonomous differential
equation
\begin{equation}
f^{\prime\prime\prime}+\alpha  ff^{\prime\prime}-\beta  f^{\prime2}=0 \label{equation}%
\end{equation}
with the boundary conditions%
\begin{equation}
f(0)=-\gamma,\ f^{\prime}\left(  \infty\right)  =0 \text{ and } f'(0)=1, \label{cond01a}%
\end{equation}%
or
\begin{equation}
f(0)=-\gamma,\ f^{\prime}\left(  \infty\right)  =0 \text{ and } f^{\prime\prime}(0)=-1. \label{cond01b}%
\end{equation}%
The first set of boundary conditions (\ref{cond01a}) with $\alpha=\frac{m+1}{2}$ and $\beta=m$ for $m \in \mathbb{R}$ corresponds to prescribed heat on the plate as in
\cite{brighicr}, \cite{brighi02}, \cite{Pop1}, \cite{cheng}, \cite{ene}, \cite{ing} and \cite{merk}. 
The second set of boundary conditions (\ref{cond01b}) with $\alpha=m+2$ and $\beta=2m+1$   for $m \in \mathbb{R}$ is for the prescribed surface heat flux as done in \cite{heat_flux} and \cite{Pop}.
In both cases the solutions depend on two parameters: $m$, the power-law exponent and $\gamma$, the mass transfer parameter. For $\gamma=0$ we have an impermeable wall, $\gamma<0$ corresponds to a fluid suction, and $\gamma>0$ to a fluid injection.

\noindent Equation $(\ref{equation})$ with suitable boundary conditions also arises in other 
industrial processes such as boundary layer flow adjacent to stretching walls 
(see \cite{banks1}, \cite{banks}, \cite{crane}, \cite{gup}, \cite{Mag}) or excitation of liquid metals
 in a high-frequency magnetic field (see \cite{mof}).
\section{The case of prescribed heat}
\subsection{Derivation of the model}

We consider a vertical permeable flat plate embedded in a porous medium at the ambient temperature $T_{\infty}$ and a rectangular Cartesian co-ordinate system with the origin fixed at the leading edge of the vertical plate, the $x$-axis directed upward along the plate and the $y$-axis normal to it. If we suppose that the porous medium is homogeneous and isotropic, that all the properties of the fluid and the porous medium are constants, that the fluid is incompressible and follows the Darcy-Boussinesq law and that the temperature along the plate is varying as $x^m$ the governing equations are 
$$\frac{\partial u}{\partial x}+\frac{\partial v}{\partial y}=0,$$
$$u=-\frac{k}{\mu}\left(\frac{\partial p}{\partial x}+\rho g \right),$$
$$v=-\frac{k}{\mu}\frac{\partial p}{\partial y},$$
$$u\frac{\partial T}{\partial x}+v\frac{\partial T}{\partial y}=
\lambda \left(\frac{\partial^2 T}{\partial x^2}
+\frac{\partial^2 T}{\partial y^2}\right ),$$
$$\rho=\rho_{\infty}(1-\beta(T-T_{\infty}))$$
where $u$ and $v$ are the Darcy velocities in the $x$ and $y$ directions, $\rho$, $\mu$ and
$\beta$ are the density, viscosity and thermal expansion coefficient of the fluid, $k$ is the
permeability of the saturated porous medium and $\lambda$ its thermal diffusivity, $p$
is the pressure, $T$ the temperature and $g$ the acceleration of the gravity. The subscript 
$\infty$ is used for a value taken far from the plate. In our system of co-ordinates the boundary conditions along the plate are
$$v(x,0)=\omega x^\frac{m-1}{2}, \quad T(x,0)=T_w(x)=T_\infty+Ax^m,\quad m\in\mathbb{R},$$
with $A>0$ and $\omega \in \mathbb{R}$ ($\omega<0$ corresponds to a fluid suction, $\omega=0$ is for an impermeable wall and $\omega>0$ corresponds to a fluid injection). The boundary conditions far from the plate are
$$u(x,\infty)=0,\quad T(x,\infty)=T_{\infty}.$$
If we introduce the stream function $\psi$ such that
$$u=\frac{\partial \psi}{\partial y}, \quad v=-\frac{\partial \psi}{\partial x}$$
and assuming that convection takes place in a thin layer around the heating plate, we obtain
the boundary layer approximation
\begin{equation}
\frac{\partial^2 \psi}{\partial y^2}
=\frac{\rho_\infty \beta g k}{\mu}\frac{\partial T}{\partial y}, \label{psi}
\end{equation}
\begin{equation}
\frac{\partial^2 T}{\partial y^2}=\frac{1}{\lambda}\left (
\frac{\partial T}{\partial x}
\frac{\partial \psi}{\partial y}-\frac{\partial T}{\partial y}
\frac{\partial \psi}{\partial x} \right ) \label{T}
\end{equation}
with 
$$\frac{\partial \psi}{\partial x}(x,0)=-\omega x^\frac{m-1}{2} \text{ and } \frac{\partial \psi}{\partial y}(x,\infty)=0.$$
Let us introduce the new dimensionless similarity variables
$$t=(Ra_{x})^\frac{1}{2}\frac{y}{x},\quad
\psi(x,y)=\lambda (Ra_{x})^\frac{1}{2}f(t),\quad
\theta(t)=\frac{T(x,y)-T_\infty}{T_w(x)-T_\infty}$$
with $Ra_x=(\rho_\infty \beta g k(T_w(x)-T_\infty)x)/(\mu\lambda)$
the local Rayleigh number. In terms of these variables equations (\ref{psi}) and (\ref{T}) become
\begin{equation}
f''-\theta'=0,
\label{i1}
\end{equation}
and
\begin{equation*}
\theta''+\frac{m+1}{2}f\theta'-mf'\theta=0,
\end{equation*}
with the boundary conditions
$$f(0)=-\gamma,\quad \theta(0)=1,$$
and
\begin{equation}
f'(\infty)=0,\quad \theta(\infty)=0,
\label{c1}
\end{equation}
where the prime denotes differentiation with respect to $t$ and
$$\gamma=\frac{2\omega}{m+1}\sqrt{\frac{\mu}{\rho_\infty \beta g k A \lambda}}.$$
Integrating (\ref{i1}) and taking into
account the boundary conditions (\ref{c1}) leads to
$$f'=\theta$$
and the problem (\ref{equation})-(\ref{cond01a}) with  $\alpha=\frac{m+1}{2}$ and $\beta=m$ for $m \in \mathbb{R}$ follows.

\subsection{Usefull tools}
\subsubsection{The initial value problem}

Let $P_{m,\gamma,\alpha}$ be the following initial value problem
\begin{equation}%
\left \{ \begin{array}
[c]{rll}%
&f^{\prime\prime\prime}+\frac{m+1}{2} ff^{\prime\prime}-m  f^{\prime2} = 0\\
&f(0) = -\gamma,\\
&f^{\prime}\left(  0\right)=1,\\
&f^{\prime\prime}(0) = \alpha.
\end{array} \right .
\end{equation}
This first approach used is a shooting method that consists in finding 
values of $f''(0)=\alpha$ for which $f$ exists on $[0,\infty)$ and such that
$f'(\infty)=0$. This direct method allows us to consider vanishing solutions but
does fail in some cases (see \cite{brighi01}).  

\subsubsection{The blowing-up co-ordinates}

Let us notice  that if $f$ is a solution of (\ref{equation}) then for all $\kappa>0$ the function $t\longmapsto\kappa f(\kappa t)$ is a solution too. Then, considering an interval $I$ on which a solution $f$ of (\ref{equation}) does not vanish, for $\tau \in I$ we can introduce the following blowing-up co-ordinates 
\begin{equation}
\forall t\in I,~~~~~s=\int_\tau^tf(\xi)d\xi,~~~~~u(s)=\frac{f'(t)}{f(t)^2}~~~~~\text{ and
}~~~~~v(s)=\frac{f''(t)}{f(t)^3}. \label{eq5}
\end{equation}
Then, we easily get
\begin{equation}
\left\{
\begin{array}[c]{l}
\dot u=P(u,v):=v-2u^2,\\
\dot v=Q_m(u,v):=-\frac{m+1}{2}v+m u^2-3uv,
\end{array}
\right.
\label{eq6}
\end{equation}
where the dot is for differentiating with respect to the variable $s$.
To come back to the original problem it is sufficient to consider the initial value problem
$P_{m,\gamma,\alpha}$ with $\gamma\not=0$ and look at the trajectories of the corresponding plane dynamical system $(\ref{eq6})$. For details, see \cite{BrighiSari}.

\subsection{Main results}

The  problem $(\ref{equation})$-$(\ref{cond01a})$ appears in engineering and physical
litterature, in very different context, in the middle of the previous century.

\noindent Rigorous mathematical results arise around the sixties. In \cite{stu}
(Appendix 2) it is mentionned that for $\gamma=0$ a simple explicit
solution can be obtained in both cases $m=1$ and $m=-\frac{1}{3}$ (see
also \cite{crane}, \cite{ing}, \cite{brighicr} and \cite{brighi02}). On the other hand, the author notes that Mr J. Watson has given a simple proof  that $(\ref{equation})$-$(\ref{cond01a})$ has no solution for $\gamma=0$ and $m\leq -1$.

\noindent An explicit solution is also given for $m=1$ and any $\gamma$, first
in \cite{gup}, and later in \cite{Mag} and \cite{brighi01}.
In these latter papers one can also find the explicit solution for
$m=-{1\over 3}$ and any $\gamma$.

\noindent Nonexistence for $\gamma=0$ and $m=-\frac{1}{2}$ was noted in \cite{banks1}. In
\cite{ing}, it is shown that for $\gamma=0$ and $m<-\frac{1}{2}$ there are no
solutions satisfying $f'f^2 \to 0$ at infinity.

\noindent Recently, further mathematical results concerning existence,
nonexistence, uniqueness, nonuniqueness and asymptotic behaviour, are
obtained in \cite{brighicr}, \cite{brighi02} for $\gamma=0$, and in \cite{guedda}, \cite{brighi01}, \cite{BrighiSari}, \cite{guedda1} and \cite{BBequiv} for the general case.

\noindent Numerical  investigations can be found in \cite{banks1}, \cite{brighi04}, \cite{Pop1}, \cite{cheng}, \cite{ing}, \cite{Mag} and \cite{wood}.

\noindent In view of all these papers, the following conclusions can be drawn.

\begin{itemize}  
\item For $m<-1$, there exists $\gamma_*>0$ such that  problem $(\ref{equation})$-$(\ref{cond01a})$ has infinitely many solutions if $\gamma>\gamma_*$, one and only
one solution if $\gamma=\gamma_*$, and no solution if $\gamma<\gamma_*$. 
For $\gamma=\gamma_{*}$ we have that $f(t)\rightarrow \lambda <0$
as $t\rightarrow \infty$, and for every $\gamma> \gamma_{*}$ there are two solutions $f$
such that $f(t)\rightarrow \lambda <0$ as $t\rightarrow \infty$ and all the other solutions verify 
$f(t)\rightarrow 0$ as $t\rightarrow \infty$
Moreover, if $f$ is a solution to $(\ref{equation})$-$(\ref{cond01a})$, then $f$ is negative, strictly increasing and either concave or convex-concave. (See Fig. 1 for the two solutions such that $f(t)\rightarrow \lambda <0$ as $t\rightarrow \infty$ and three other solutions in the case $m=-2$ and $\gamma=5$.)
$$\begin{array}
[c]{cc}
\hspace{0.2cm} & \includegraphics[scale=0.7,trim=0 480 0 0]{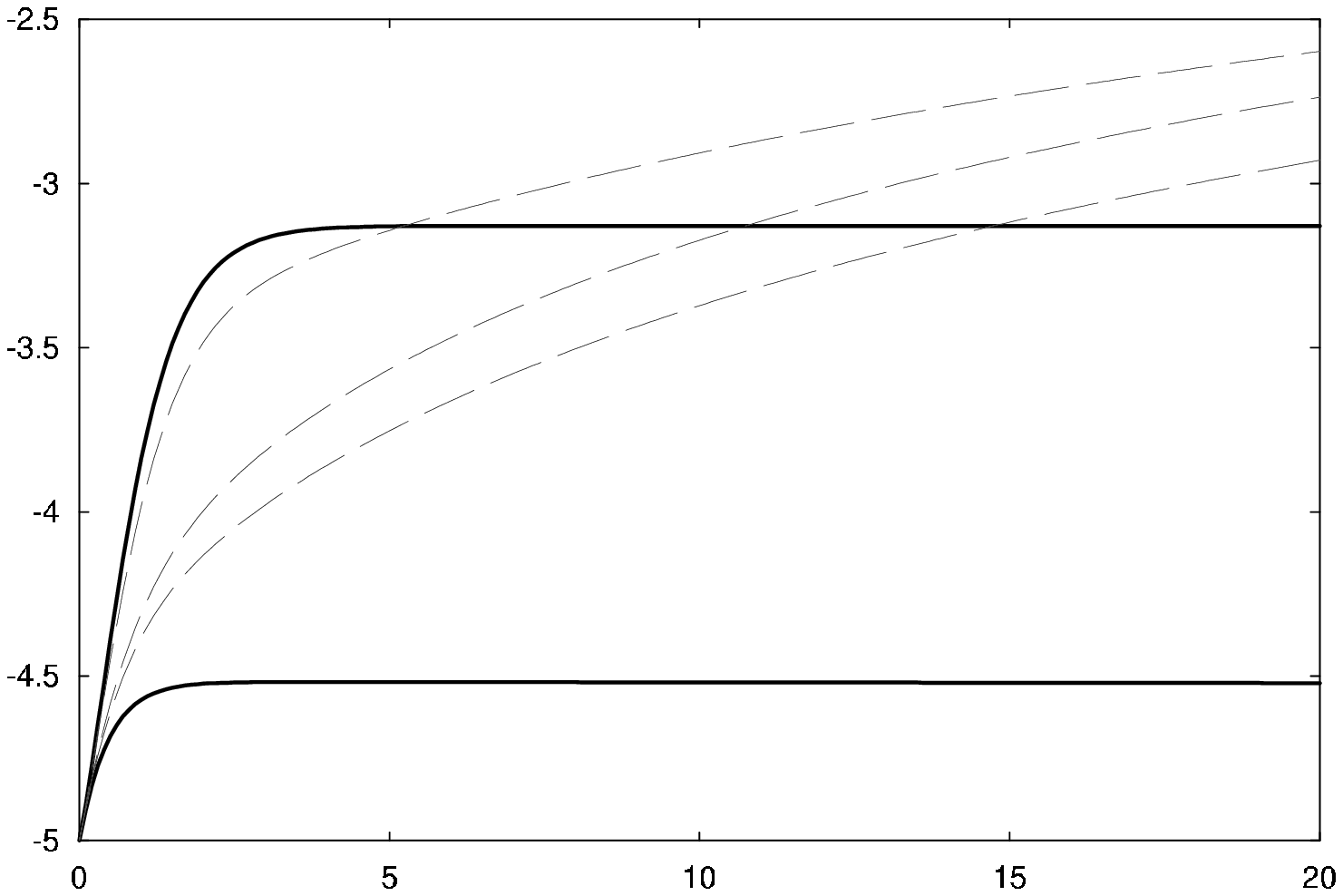}
\end{array}$$
\vspace{-19mm}

\centerline{Fig. 1}   
\item For $m=-1$ and for every $\gamma\in\mathbb{R}$, the problem  $(\ref{equation})$-$(\ref{cond01a})$ has no solution.
\item For $-1<m\leq-{1\over 2}$ and for every $\gamma\geq 0$, the
problem  $(\ref{equation})$-$(\ref{cond01a})$ has no solution.
\item For $-1<m<-{1\over 2}$, there exists $\gamma_*<0$ such that problem
 $(\ref{equation})$-$(\ref{cond01a})$ has no solution for $\gamma_*<\gamma<0$, one and only one solution which is bounded for $\gamma=\gamma_*$, and two bounded solutions and infinitely many unbounded solutions for $\gamma<\gamma_*$. These solutions are strictly increasing and either concave or convex-concave. (See Fig. 2 to see the two bounded solutions and four unbouded solutions in the case $m=-0.75$ and $\gamma=-10$.)

\item For $-{1\over 2}\leq m<-{1\over 3}$ and for every $\gamma<0$, the problem
 $(\ref{equation})$-$(\ref{cond01a})$ has one bounded solution and infinitely many unbounded solutions. All these solutions are strictly increasing and either concave or convex-concave.
\item For $-{1\over 3}\leq m<0$ and for every $\gamma\in\mathbb{R}$, the problem
 $(\ref{equation})$-$(\ref{cond01a})$ has an infinite number of solutions. Moreover, if $\gamma\leq 0$ one and only one solution is bounded, and if $\gamma>0$ at least one is bounded, many infinitely are unbounded. All solutions are strictly increasing and either concave or convex-concave. If $\gamma>0$, the solutions becomes positive for large $t$.
 $$\begin{array}
[c]{cc}
\hspace{0.2cm} & \includegraphics[scale=0.7,trim=0 480 0 0]{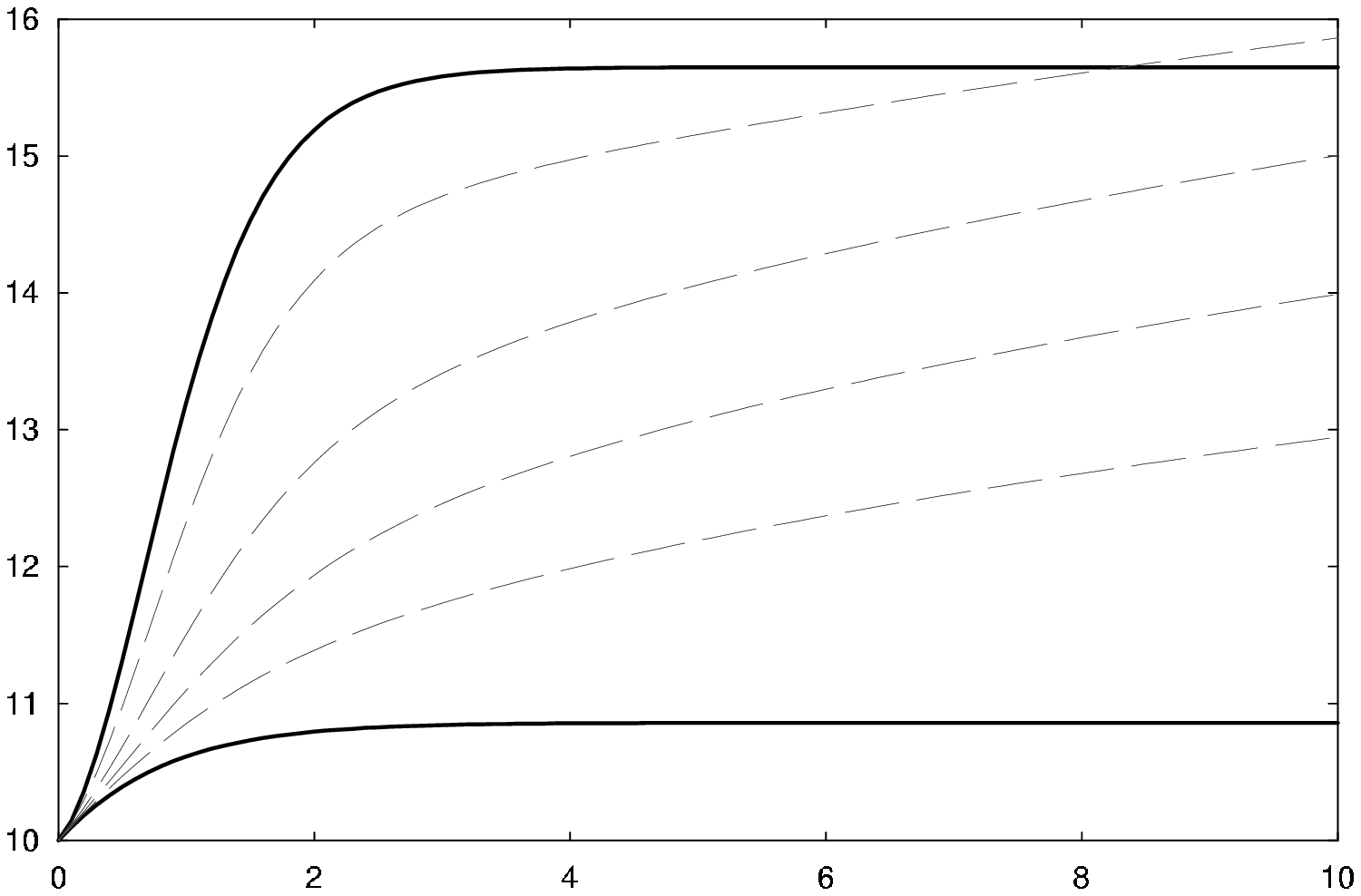}
\end{array}$$
\vspace{-19mm}

\centerline{Fig. 2} 
\item For $ m\in[0,1]$ and for every $\gamma\in\mathbb{R}$, the problem  $(\ref{equation})$-$(\ref{cond01a})$ has one and only one solution, moreover this solution is concave and 
bounded. (See Fig. 3 for the unique solution in the case $m=0.5$ and $\gamma=0$.)
$$\begin{array}
[c]{cc}
\hspace{0.2cm} & \includegraphics[scale=0.7,trim=0 480 0 0]{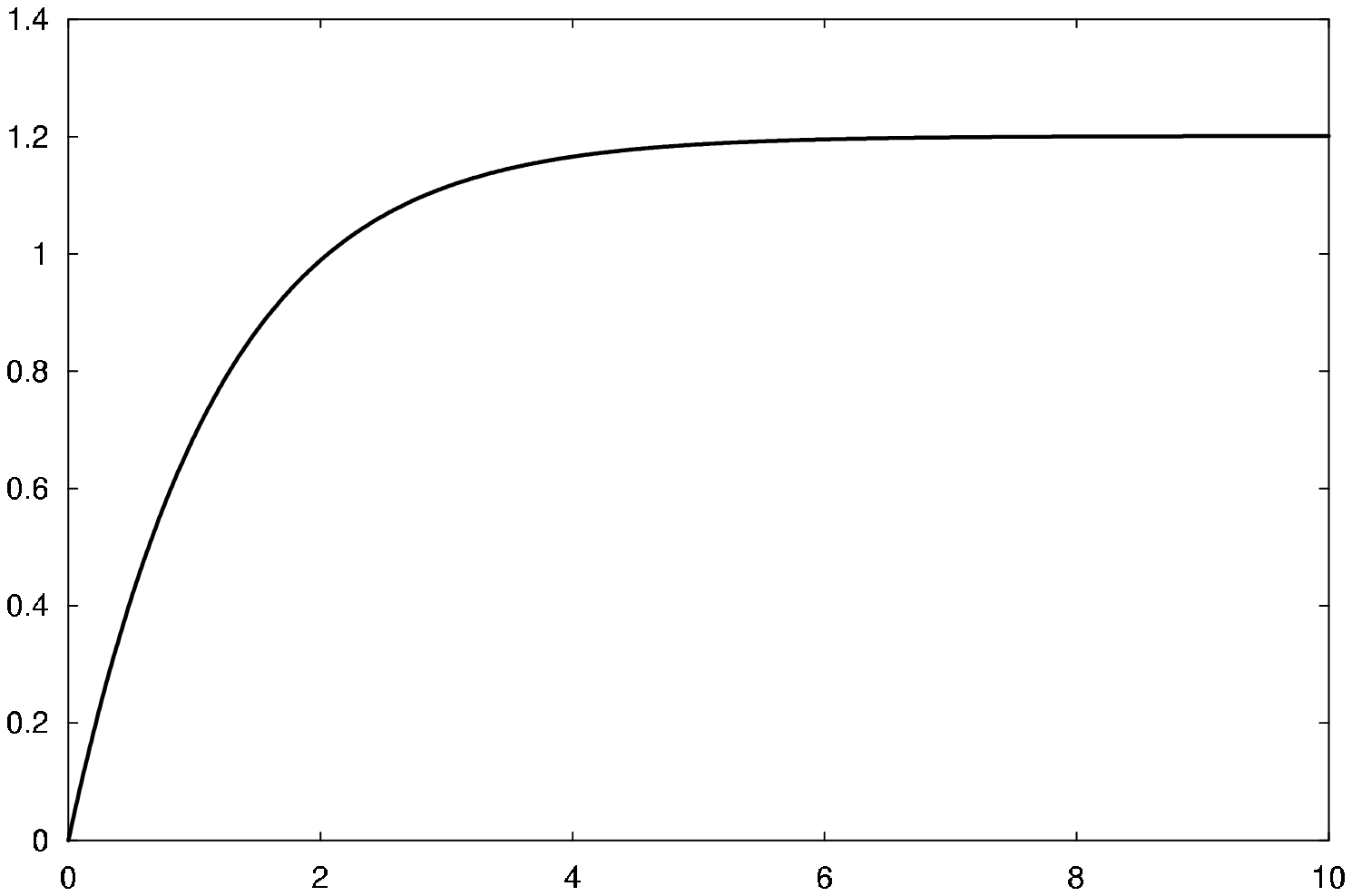}
\end{array}$$
\vspace{-19mm}

\centerline{Fig. 3}   
\item For $ m>1$ and for every $\gamma\in \mathbb{R}$, the problem  $(\ref{equation})$-$(\ref{cond01a})$ has one and only one concave solution and an infinite number of concave-convex solutions. All these solutions are bounded. Moreover, there is an unique concave-convex solution that verifies $f(t)\rightarrow \lambda>0$ as $t\rightarrow \infty$ and all the other concave-convex solutions are such that $f(t)\rightarrow 0$ as $t\rightarrow \infty$.
(See Fig. 4 for the unique concave solution and three concave-convex solutions in the case $m=1.1$ and $\gamma=0$.)
$$\begin{array}
[c]{cc}
\hspace{0.2cm} & \includegraphics[scale=0.7,trim=0 480 0 0]{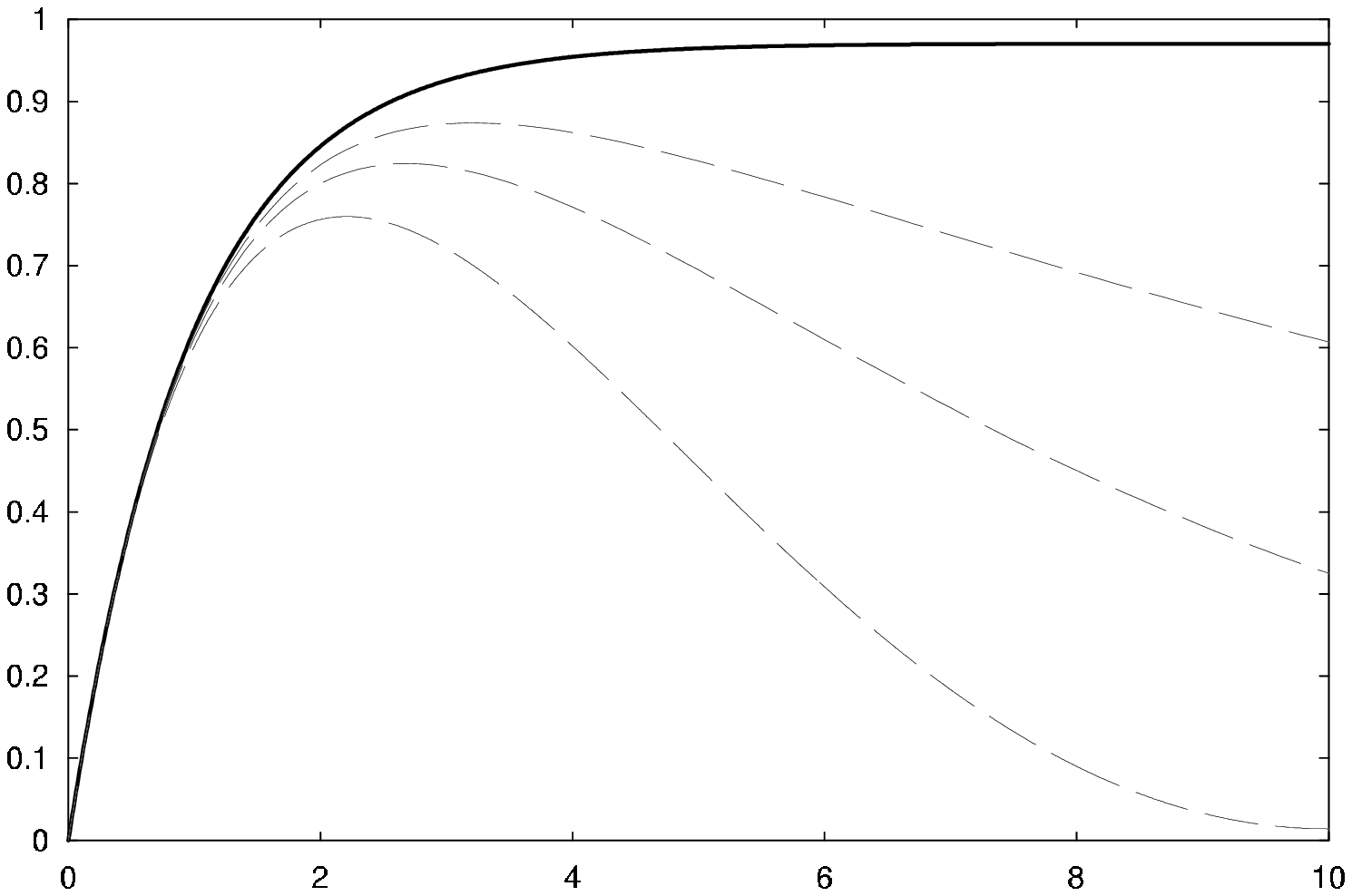}
\end{array}$$
\vspace{-19mm}

\centerline{Fig. 4}   
\end{itemize}

\begin{remark}
The case $m=0$ leads to the well know Blasius equation (see {\rm \cite{brighi03}}, 
{\rm \cite{bla}}) that also is a special case of the Falkner-Skan equation (see {\rm \cite{falk}}).
\end{remark}
\begin{remark}
In {\rm \cite{Mag1}} the authors gives results about a slightly different problem for $m=-1$ that involves pseudo-similarity.
\end{remark}

\noindent We see from these results, that the unsolved questions concern the case $\gamma\geq 0$. More precisely, it should be interesting to try to answer to the following points
\begin{itemize}
\item For $-{1\over 2}< m<-{1\over 3}$, what happens for $\gamma\geq 0$ ?
\item For $-{1\over 3}\leq m<0$ and $\gamma>0$, is there one or more bounded solutions ?
\end{itemize}
Another purpose is to compute the critical values $\gamma_*$ appearing in the results above. 

\section{The case of prescribed heat flux}
We now suppose that the plate is subjected to a variable heat flux varying as $x^m$ and a mass transfer rate varying as $x^\frac{m-1}{3}$ following \cite{Pop} to obtain the problem 
(\ref{equation})-(\ref{cond01b}) with $\alpha=m+2$ and $\beta=2m+1$. The mathematical study is made in \cite{heat_flux} and leads to the following results

\begin{itemize}

\item For $m<-2$ there exists 
$\gamma_{*}>\sqrt[3]{\frac{2}{(m+2)^2}}$ such that the problem
$\left(  \ref{equation}\right)  $-$\left(  \ref{cond01b}\right)  $ has
no solution for $\gamma<\gamma_{\ast}$, one and only one solution
for $\gamma=\gamma_{\ast}$ and infinitely many solutions for $\gamma>\gamma_{\ast}.$
For $\gamma=\gamma_{*}$ we have that $f(t)\rightarrow \lambda <0$
as $t\rightarrow \infty$, and for every $\gamma> \gamma_{*}$ there are two solutions $f$
such that $f(t)\rightarrow \lambda <0$ as $t\rightarrow \infty$ and all the other solutions verify 
$f(t)\rightarrow 0$ as $t\rightarrow \infty$
Moreover, if $f$ is a solution of $\left(  \ref{equation}\right)  $-$\left(
\ref{cond01b}\right)  $, then $f$ is negative, strictly concave and
increasing.

\item For $m=-2$ and for every $\gamma\in\mathbb{R}$, the problem $\left(
\ref{equation}\right)  $-$\left(  \ref{cond01b}\right)  $ has no solution.

\item For $-2<m<-1$, there exists $\gamma_{\ast}<0$ such that the problem
$\left(  \ref{equation}\right)  $-$\left(  \ref{cond01b}\right)  $ has
no solution for $\gamma>\gamma_{\ast}$, one and only one solution which is bounded 
for $\gamma=\gamma_{\ast}$ and two bounded solutions and infinitely many 
unbounded solutions for $\gamma<\gamma_{\ast}.$
Moreover, if $f$ is a solution of $\left(  \ref{equation}\right)  $-$\left(
\ref{cond01b}\right) $, then $f$ is positive, strictly concave,
increasing and $f^{\prime}(0)\geq -\frac{1}{(m+2)\gamma}.$

\item For $m=-1$ the problem $\left(  \ref{equation}\right)  $-$\left(
\ref{cond01b}\right) $ only admits solutions for $\gamma<0$. In this case there is an unique 
bounded solution with $f'(0)=-\frac{1}{\gamma}$ and an infinite 
number of unbounded solutions with $f'(0)>-\frac{1}{\gamma}$. Moreover all the
solutions are positive, strictly concave and increasing.

\item For $-1< m<-\frac{1}{2}$ the problem $\left(  \ref{equation}\right)  $-$\left(
\ref{cond01b}\right) $ admits at least one bounded solution for $\gamma \in \mathbb{R}$ and 
many infinitely unbounded solutions for $\gamma<0$. All these solutions are increasing and 
strictly concave and uniqueness of the bounded solution holds for $\gamma\leq 0$.

\item For $m\geq -\frac{1}{2}$ all the solutions are bounded.

\item For $-\frac{1}{2}\leq m \leq 1$ and for every $\gamma \in \mathbb{R}$ the problem 
(\ref{equation})-(\ref{cond01b}) has one and only one solution. This solution is strictly concave
and increasing. (Let us notice that for $m=-\frac{1}{2}$ we have the Blasius equation.)

\item For $m>1$ and $\gamma \in \mathbb{R}$ the problem (\ref{equation})-(\ref{cond01b}) 
has one and only one concave solution and infinitely many concave-convex solutions. 
Moreover, there is an unique concave-convex solution 
that verifies $f(t)\rightarrow \lambda>0$ as $t\rightarrow \infty$ and all the other 
concave-convex solutions are such that $f(t)\rightarrow 0$ as $t\rightarrow \infty$.
\bigskip
\end{itemize}

In this case it remains only two open questions
\begin{itemize}
\item For $-1<m< -\frac{1}{2}$ and $\gamma >0$, is the bounded solution unique ?
\item For $-1<m< -\frac{1}{2}$ and $\gamma \geq 0$, is there unbounded solution ?
\end{itemize}

\section{Asymptotic behaviour of the unbounded solutions}
For the equation (\ref{equation}) we have the following asymptotic equivalent found in 
\cite{BBequiv} and \cite{guedda1} that holds for unbounded solutions
\begin{theorem} Let $f$ be an unbounded solution of
$(\ref{equation})$-$(\ref{cond01a})$ or $(\ref{equation})$-$(\ref{cond01b})$. There exists a constant $c>0$ such that  $$|f(t)| \sim c t^\frac{\alpha}{\alpha-\beta}\quad \text{ as }\quad t\to\infty.$$
\end{theorem}

\end{document}